\documentclass[11pt]{amsart}
\usepackage{amssymb,amsmath,amsgen,amsxtra,amsfonts,a4} 
\usepackage{dsfont}
\usepackage[matrix,arrow]{xy}

\begin{document}
%
%
%
\theoremstyle{definition}
\newtheorem{Definition}{Definition}[section]
\newtheorem{Convention}{Definition}[section]
\newtheorem{Construction}{Construction}[section]
\newtheorem{Example}[Definition]{Example}
\newtheorem{Examples}[Definition]{Examples}
\newtheorem{Remark}[Definition]{Remark}
\newtheorem{Draft}[Definition]{Remark for draft}
\newtheorem{Remarks}[Definition]{Remarks}
\newtheorem{Caution}[Definition]{Caution}
\newtheorem{Conjecture}[Definition]{Conjecture}
\newtheorem*{Question}{Question}
\newtheorem*{Answer}{Answer}
\newtheorem*{Acknowledgement}{Acknowledgement}
\theoremstyle{plain}
\newtheorem{Theorem}[Definition]{Theorem}
\newtheorem*{Theoremx}{Theorem}
\newtheorem{Proposition}[Definition]{Proposition}
\newtheorem{Lemma}[Definition]{Lemma}
\newtheorem{Corollary}[Definition]{Corollary}
\newtheorem{Fact}[Definition]{Fact}
\newtheorem{Facts}[Definition]{Facts}
\newtheoremstyle{voiditstyle}{3pt}{3pt}{\itshape}{\parindent}%
{\bfseries}{.}{ }{\thmnote{#3}}%
\theoremstyle{voiditstyle}
\newtheorem*{VoidItalic}{}
\newtheoremstyle{voidromstyle}{3pt}{3pt}{\rm}{\parindent}%
{\bfseries}{.}{ }{\thmnote{#3}}%
\theoremstyle{voidromstyle}
\newtheorem*{VoidRoman}{}
%
\newcommand{\prf}{\par\noindent{\sc Proof.}\quad}
\newcommand{\blowup}{\rule[-3mm]{0mm}{0mm}}
\newcommand{\cal}{\mathcal}
\newcommand{\Aff}{{\mathds{A}}}
\newcommand{\BB}{{\mathds{B}}}
\newcommand{\CC}{{\mathds{C}}}
\newcommand{\FF}{{\mathds{F}}}
\newcommand{\GG}{{\mathds{G}}}
\newcommand{\HH}{{\mathds{H}}}
\newcommand{\NN}{{\mathds{N}}}
\newcommand{\ZZ}{{\mathds{Z}}}
\newcommand{\PP}{{\mathds{P}}}
\newcommand{\QQ}{{\mathds{Q}}}
\newcommand{\RR}{{\mathds{R}}}
\newcommand{\Liea}{{\mathfrak a}}
\newcommand{\Lieb}{{\mathfrak b}}
\newcommand{\Lieg}{{\mathfrak g}}
\newcommand{\Liem}{{\mathfrak m}}
\newcommand{\ideala}{{\mathfrak a}}
\newcommand{\idealb}{{\mathfrak b}}
\newcommand{\idealg}{{\mathfrak g}}
\newcommand{\idealm}{{\mathfrak m}}
\newcommand{\idealp}{{\mathfrak p}}
\newcommand{\idealq}{{\mathfrak q}}
\newcommand{\idealI}{{\cal I}}
\newcommand{\lin}{\sim}
\newcommand{\num}{\equiv}
\newcommand{\dual}{\ast}
\newcommand{\iso}{\cong}
\newcommand{\homeo}{\approx}
\newcommand{\mm}{{\mathfrak m}}
\newcommand{\pp}{{\mathfrak p}}
\newcommand{\qq}{{\mathfrak q}}
\newcommand{\rr}{{\mathfrak r}}
\newcommand{\pP}{{\mathfrak P}}
\newcommand{\qQ}{{\mathfrak Q}}
\newcommand{\rR}{{\mathfrak R}}
%
%
\newcommand{\dq}{{``}}
\newcommand{\OO}{{\cal O}}
\newcommand{\into}{{\hookrightarrow}}
\newcommand{\onto}{{\twoheadrightarrow}}
\newcommand{\Spec}{{\rm Spec}\:}
\newcommand{\Proj}{{\rm Proj}\:}
\newcommand{\Pic}{{\rm Pic }}
\newcommand{\Br}{{\rm Br}}
\newcommand{\NS}{{\rm NS}}
\newcommand{\chit}{\chi_{\rm top}}
\newcommand{\KanDiv}{{\cal K}}
\newcommand{\perdef}{{\stackrel{{\rm def}}{=}}}
\newcommand{\Cycl}[1]{{\ZZ/{#1}\ZZ}}
\newcommand{\Sym}{{\mathfrak S}}
\newcommand{\Xcan}{X_{{\rm can}}}
\newcommand{\ab}{{\rm ab}}
\newcommand{\Aut}{{\rm Aut}}
\newcommand{\Hom}{{\rm Hom}}
\newcommand{\Supp}{{\rm Supp}}
\newcommand{\ord}{{\rm ord}}
\newcommand{\divisor}{{\rm div}}
\newcommand{\Alb}{{\rm Alb}}
\newcommand{\Jac}{{\rm Jac}}
\newcommand{\piet}{{\pi_1^{\rm \acute{e}t}}}
\newcommand{\Het}[1]{{H_{\rm \acute{e}t}^{{#1}}}}
\newcommand{\Hcris}[1]{{H_{\rm cris}^{{#1}}}}
\newcommand{\HdR}[1]{{H_{\rm dR}^{{#1}}}}
\newcommand{\hdR}[1]{{h_{\rm dR}^{{#1}}}}
\newcommand{\defin}[1]{{\bf #1}}
\title[Birational modifications of surfaces]{Birational modifications of surfaces via unprojections}
\author{Christian Liedtke\and Stavros Argyrios Papadakis}
\dedicatory{March 1, 2010}

\begin{abstract}
We describe elementary transformations between
minimal models of rational surfaces in terms of
unprojections.
These do not fit into the 
framework of Kustin--Miller unprojections as introduced
by Papadakis and Reid, since we have to leave the world
of projectively Gorenstein varieties.
Also, our unprojections do not depend on the choice
of the unprojection locus only, but need extra data 
corresponding to the choice of a divisor on this unprojection
locus.
\end{abstract}
\setcounter{tocdepth}{1}
\maketitle
%
\section*{Introduction}

Unprojection, introduced by Miles Reid in \cite{kinosaki},
is a technique to describe birational transformations in higher dimensional
geometry explicitly in terms of commutative algebra.
As explained in \cite[Section 2.3]{pr},
the prototype and easiest example of an unprojection is the Castelnuovo
blow-down of a rational $(-1)$-curve lying on a smooth cubic surface
in $\PP^3$ as the inverse of a projection from a del~Pezzo surface
of degree $4$ in $\PP^4$, which also explains the name unprojection.

Since then, unprojections have been applied to many explicit constructions
in birational geometry of surfaces and $3$-folds,
compare \cite{abm}, \cite{cpr}, \cite{rs},...
Moreover, these techniques have also proved to be useful for the
construction of key varieties, see, for example,
\cite{kinosaki} or \cite{campedelli}.
Finally, the general theory of unprojection has been developed further
by the second author and a general framework has been proposed
in \cite{towards}.

In this note we return to the two-dimensional case.
Here, every birational map between two smooth projective
surfaces can be factored into a sequence of blow-ups of points and 
Castelnuovo blow-downs of rational $(-1)$-curves.
A prominent class of birational transformations 
between smooth surfaces 
are the {\it elementary transformations},
which relate minimal models of surfaces of Kodaira dimension
$\kappa=-\infty$ to each other, compare \cite[Chapter 12]{bad}.
Their higher-dimensional generalisations are Sarkisov links,
compare \cite{cpr}.

\begin{Question}
  Can one describe elementary transformations of surfaces in 
  terms of unprojections?
\end{Question}

For reasons of simplicity and since all relevant problems
already occur within this class of surfaces, we will
restrict ourselves to minimal rational surfaces.
Already here we encounter new types of unprojections
and new phenomena show up.
This is related to the fact that an elementary transformation 
depends not only on the choice of a curve but needs also 
the choice of a point on it.

\begin{Answer}
   Yes, this can be done for minimal rational surfaces.
   However, the unprojection is no longer determined by
   the unprojection locus alone.
   Moreover, this cannot be done within the framework of
   projectively Gorenstein varieties, as in the classical
   case of Kustin--Miller unprojections.
\end{Answer}

Minimal rational surfaces consist of $\PP^2$ and 
Hirzebruch surfaces.
By definition, the Hirzebruch surface, sometimes also
called Segre surface, $\FF_d$ 
is the $\PP^1$-bundle 
$\PP(\OO_{\PP^1}\oplus\OO_{\PP^1}(d))\to\PP^1$.
An elementary transformation of $\FF_d$ is the following:
we choose a point lying on a fibre of this projection
and blow it up.
The strict transform of this fibre on the blow-up is
a rational $(-1)$-curve and blowing it down we obtain
the desired elementary transformation of $\FF_d$.
Depending on the position of the point we blew up to
start with, the resulting surface is isomorphic to 
$\FF_{d+1}$ or $\FF_{d-1}$.

As embedding for the Hirzebruch surfaces we choose
their realisations as surfaces of minimal degree, i.e.,
as scrolls.
As unprojection locus $\Gamma$ we choose a line lying on this
scroll, which corresponds to the fibre of the projection
of this Hirzebruch surface onto $\PP^1$.
In terms of rings and ideals we have
$$
V(I)\,=\,\Gamma\,\subset\,\Proj S\,\subset\,
\PP^{N}\,.
$$
The point of departure for unprojections is 
the $S$-module $\Hom_S(I,S)$, which in our  
situation turns out to
be generated by two elements in degree zero.
This implies that the associated unprojection ring 
is ''not geometric``, although somewhat similar 
rings have been considered in \cite{kinosaki} 
to describe $3$-fold flips.
Geometrically, this is related to the fact that
unprojections correspond to contractions and
that $\Gamma$ has not negative self-intersection,
which would be necessary in order to contract it.

Instead, we will use natural submodules of 
$\Hom_S(I,S)$ to construct our unprojections.
Geometrically, these submodules correspond to 
choosing a divisor $D$ on $\Gamma=V(I)$.
In case $D$ is a divisor of degree $k\geq1$, whose
support consists of $k$ {\it distinct} points,
our results specialise to the following

\begin{Theoremx}
  The unprojection of 
  $X=\Proj S\subset\PP^N$ 
  with respect to $D\subset\Gamma\subset X$ 
  is a normal and projectively Cohen--Macaulay
  surface inside $\PP(1^{N+1},k)$,
  which arises from $X$ by first blowing up $D$ and then
  contracting the strict transform of
  $\Gamma$. 
\end{Theoremx}

In this special case the unprojection is smooth outside a
toric singularity of type $\frac{1}{k}(1,1)$, which is
induced from the singularity of the ambient weighted 
projective space.
We refer to Section \ref{sec:unprojection} for the general case.

The case $k=1$ corresponds to elementary transformations
of Hirzebruch surfaces.
Moreover, if $\Proj S\iso\FF_d$ and we vary the divisor $D$,
which is just one point in this case, we obtain 
a $1$-parameter family of unprojections, all of which
are isomorphic to $\FF_{d-1}$ except one surface which is
isomorphic to $\FF_{d+1}$.
This fits nicely into the deformation and degeneration 
theory of Hirzebruch surfaces, confer 
\cite[Theorem VI.8.1]{bhpv}.

Finally, we give an application to odd Horikawa surfaces,
i.e., to  minimal surfaces of general type with 
$K^2=2p_g-3$.
More precisely, for an odd Horikawa surface with $p_g\geq7$
the canonical and the bicanonical image are rational
surfaces.
Moreover, the canonical image is a Hirzebruch surface
realised as surface of minimal degree.
Then the birational transformation relating canonical
and bicanonical image corresponds to an 
unprojection of the type considered in this article
with $k=2$.

\begin{Acknowledgement}
  We thank Francesco~Zucconi and the referee for remarks, 
  comments and pointing out a couple of inaccuracies.
  Stavros~Papadakis is a participant of  the  Project
  PTDC/MAT/099275/2008, a member of  CAMSGD (IST/UTL),
  and gratefully acknowledges funding from the
  Portuguese Funda\c{c}\~ao para a Ci\^encia e a Tecno\-lo\-gia
  (FCT) under  research grant SFRH/BPD/22846/2005 of POCI2010/FEDER. 
  Much of this work was done while he was visiting 
  the university of D\"usseldorf. 
  We thank the Mathe\-mati\-sches Institut and the DFG-Forschergruppe
  {\it Classification of Algebraic Surfaces and Compact Complex Mani\-folds}
  for the kind hospitality and financial support.
\end{Acknowledgement}

\section{Hirzebruch surfaces and scrolls}
\label{sec:scrolls}

We fix once and for all an arbitrary field $k$ over which 
all our schemes will be defined.
Let $d\geq0$ be a non-negative integer.
Then the {\it Hirzebruch surface}, or,
{\it Segre surface}, $\FF_d$ is defined to be the
$\PP^1$-bundle $\PP(\OO_{\PP^1}\oplus\OO_{\PP^1}(d))\to\PP^1$.
We denote by $\Gamma$ the class of a fibre of this projection.
Moreover, there exists a section $\Delta_0$ with self-intersection
$-d$, which is unique if $d\neq0$.

We remark that $\FF_0$ is isomorphic to $\PP^1\times\PP^1$,
and that $\FF_1$ is isomorphic to $\PP^2$ blown-up in 
one point. 
Moreover, $\FF_1$ is the only Hirzebruch that is not 
minimal, and among Hirzebruch surfaces, the only del~Pezzo
surfaces are $\FF_0$ and $\FF_1$. 
Now, a projectively Cohen--Macaulay scheme $X$ is 
projectively Gorenstein if and only if there exists a 
$k\in\ZZ$ such that $\omega_X\iso\OO_X(k)$,
see \cite[Section 21.11]{eis}.
Thus, apart from $\FF_0$ and $\FF_1$, Hirzebruch
surfaces do {\it not} possess embeddings into 
projective space that are projectively Gorenstein.
In particular, elementary transformations
of minimal rational surfaces in terms of unprojections 
cannot be described within the framework of 
Kustin--Miller unprojections as in \cite{pr}.

However, Hirzebruch surfaces do possess nice
embeddings into projective space.
Namely,
for integers $m,n$ satisfying $n\geq m\geq1$
we define $\FF(m,n)$ to be the {\it surface scroll}
in $\PP^{m+n+1}$
defined by the vanishing of the $2\times2$-minors of
\begin{equation}
\label{scroll equation}
\left(
\begin{array}{ccc|ccc}
  x_{00} & ... & x_{0m-1} & x_{10} & ... & x_{1n-1}\\
  x_{01} & ... & x_{0m} & x_{11} & ... & x_{1n}
\end{array}
\right)\,.
\end{equation}
Abstractly, this scroll is isomorphic to $\FF_{n-m}$.
Moreover, $\FF(m,n)$ corresponds to embedding 
$\FF_{d}$ with $d=n-m$
via the complete linear system $|\Delta_0 + n\Gamma|$ 
into projective space.
Under this embedding, the projection onto $\PP^1$
is given by the ratios of the columns
of (\ref{scroll equation})
$$
\begin{array}{ccc}
 \FF(m,n) &\to& \PP^1 \\%
 \left[ x_{00}:...:x_{0m}:x_{10}:...:x_{1,n} \right]
 &\mapsto&
 \left[x_{00}:x_{01}\right]\,=\,...\,=\,\left[x_{1n-1}:x_{1n}\right]\,.
\end{array}
$$
Let us fix the fibre $\Gamma$ over $[0:1]$, which is given
by the vanishing of the first row 
in (\ref{scroll equation}):
\begin{equation}
\label{gamma equation}
\Gamma\,=\,\{ x_{00}\,=\,...\,=\,x_{0m-1}\,=\,x_{10}
\,=\,...\,=\,x_{1n-1}\,=\,0\}
\,\cap\,\FF(m,n)\,.
\end{equation}

Apart from the nice determinantal description there
is another reason why these embeddings of the Hirzebruch
surfaces are distinguished.
Namely, a non-degenerate and integral
surface in $\PP^N$ has degree at least $N-1$.
If such a surface has degree equal to $N-1$, 
then a theorem of del~Pezzo states that it is 
precisely one of the $\FF(m,n)$'s above, 
$\PP^2$, the Veronese surface in $\PP^5$,
or the cone over a rational normal curve, 
confer \cite{eh} for a modern account.
The case of the cone over a rational normal curve
corresponds to having only one block
in the matrix (\ref{scroll equation}) above.

We set $R=k[x_{0i},x_{1j}]$ with $i=0,...,m$ 
and $j=0,...,n$,
which is the homogeneous coordinate ring of $\PP^{m+n+1}$.
Let $Q$ be the ideal of $R$ corresponding to $\FF(m,n)$,
i.e., the homogeneous ideal generated by the
$2\times2$ minors of (\ref{scroll equation}).
We set $S=R/Q$ and define $I$ to be the ideal of $S$ 
defining $\Gamma$, i.e., the ideal corresponding
to (\ref{gamma equation})
$$
I\,=\,(x_{00},\,...,\,x_{0m-1},\,x_{10},\,...,\,x_{1n-1})
\,\subset\, S\,=\,R/Q\,.
$$
The ring $S$ is Cohen--Macaulay, which is related
to the fact that the embedding of $\FF(m,n)$ into  $\PP^{m+n+1}$ 
is given by a complete linear system, 
confer \cite[Exercise 18.16]{eis}.
Alternatively, it also follows from the determinantal description of $S$
by Eagon's theorem, see \cite[Theorem 18.18]{eis}.
\medskip

Since we want to unproject from $\Gamma=V(I)$, we need
an analysis of the graded $S$-module
$\Hom_S(I,S)$, compare \cite[Section 4]{towards}.
From loc.cit. or \cite[Section 1]{pr} we recall the 
short exact sequence
\begin{equation}
 \label{poincare residue}
0\,\longrightarrow\,S\,\longrightarrow\,
\Hom_S(I,S)\,\stackrel{{\rm res}}{\longrightarrow}\,
{\rm Ext}^1_S(S/I,S)\,\longrightarrow\,0\,,
\end{equation}
where ${\rm res}$ stands for Poincar\'e residue map.
Obviously, the inclusion $\imath$ of $I$ into $S$ lies
in $\Hom_S(I,S)$.
Moreover, the map
$$
\begin{array}{cccccl}
\phi&:&x_{0i} &\mapsto& x_{0i+1}& i=0,...,m-1\\
 && x_{1j} &\mapsto& x_{1j+1}& j=0,...,n-1\\
\end{array}
$$
sending the first row of (\ref{scroll equation}) to the second
row defines a homomorphism of degree zero of graded $S$-modules 
from $I$ to $S$.
This is best seen by considering the element
$\tilde{s}=x_{01}/x_{00}=...=x_{1n}/x_{1n-1}$ in the field
of fractions $k(S)$ of the domain $S$.
Then multiplication by $\tilde{s}$ induces a homomorphism of $S$-modules 
from $I$ to $k(S)$ which yields $\phi$.

\begin{Proposition}
  \label{generators}
  The $S$-module $\Hom_S(I,S)$ is generated by the two elements
  $\imath$ and $\phi$, both of which are of degree zero.
\end{Proposition}

Since this module is crucial for the construction and analysis of
unprojections, we decided to give two proofs -- one more geometric
and one purely algebraic:
\medskip

\par\noindent{\sc First Proof.}\quad
Let $X=\Proj S$ together with its very ample invertible sheaf $\OO_X(1)$.
From the determinantal description we infer that $X$ is projectively
normal, which implies $S=\bigoplus_{n\in\ZZ}H^0(X,\OO_X(n))$.
Moreover, the sheafification of ${\rm Hom}_S(I,S)$ is $\OO_X(\Gamma)$, which implies
that there is a natural injection of graded 
$S$-modules
$$
\alpha: {\rm Hom}_S(I,S) \,\to\,\bigoplus_{n\in\ZZ} H^0(X,\OO_X(\Gamma)(n))\,=:\,M.
$$
Thus, there are no elements of negative degree in ${\rm Hom}_S(I,S)$.
In degree zero, we have $\imath$ and $\phi$ 
in ${\rm Hom}_S(I,S)$ and $h^0(X,\OO_X(\Gamma))=2$, which implies
that $\alpha$ is an isomorphism in degree zero.

Using the explicit description of global sections of invertible sheaves on scrolls
in terms of bihomogeneous polynomials as in \cite[Chapter 2]{chapters}, it follows
easily that $M$ is generated as an $S$-module in degree zero.
It follows that $\alpha$ is an isomorphism of graded $S$-modules and hence
that ${\rm Hom}_S(I,S)$ is generated by $\imath$ and $\phi$.
\qed\medskip

\par\noindent{\sc Second Proof.}\quad
Denote by $B$ the following subset of $R$
$$
\begin{array}{cclclcl}
 B &=& \{1\} &\cup& \{ x_{0i}\cdot x_{0m}^a\cdot x_{1n}^b &|& 2\leq i\leq m-1,\,a,b\geq0 \} \\
 & & &\cup& \{ x_{0m}^a\cdot x_{1n}^b &|& a,b\geq0,\,(a,b)\neq(0,0) \} \\
 & & &\cup& \{ x_{10}^a\cdot x_{1n}^b &|& a,b\geq0,\,(a,b)\neq(0,0) \} \\
 & & &\cup& \{ x_{10}^a\cdot x_{1i} \cdot x_{1n}^b &|& 1\leq i\leq n-1,\, a,b\geq0 \}\,. \\
\end{array}
$$
\begin{center}
 {\sc First claim: } The set $B$ is a basis of the $k$-vector space $R/(Q + (x_{00}, x_{01}))$. 
\end{center}
Set $Q_1=Q+(x_{00},x_{01})$. 
Using the relations $x_{0i}x_{0j}=x_{0i-1}x_{0j+1}$,
$x_{1i}x_{1j}=x_{1i-1}x_{1j+1}$ and $x_{0i}x_{1j}=x_{0i-1}x_{1j+1}$ of
$Q_1$ it is 
not difficult to see that given a monomial $w\in R$
there exists another monomial $w'\in R$ with $w-w'\in Q_1$ such that
$w'\in B$.
This shows that $B$ spans $R/Q_1$ as $k$-vector space.
To prove linear independence, we consider the $k$-algebra homomorphism
$g:R\to k[z,s,t]$ defined by 
$$
g(x_{0i})\,=zt^{m-i}s^i,\,0\leq i\leq m\mbox{ \quad and \quad }
g(x_{1j})\,=\,t^{n-j}s^j,\,0\leq j\leq n\,.
$$
Now, assume that we have an element 
$$
a\,=\,\sum_{b\in B} a_b\cdot b\in Q_1,\mbox{ where almost all }a_b=0\,.
$$
Using $Q\subseteq\ker g$, we see that
$g(a)=\sum_b a_b g(b)$ lies inside the ideal of $k[z,s,t]$ 
generated by $zt^{m-1}$.
On the other hand, none of the monomials $g(b)$, $b\in B$ is divisible
by $zt^{m-1}$.
Hence if $g(a)\neq0$ we get a contradiction, so $g(a)=0$.
Since it is clear that the set $\{g(b)\,|\,b\in B\}$ is linearly independent,
we get $a_b=0$ for all $b\in B$ and we conclude linear independence of $B$.
This proves the claim.
\begin{center}
 {\sc Second claim: } If $u\in R$ fulfills $ux_{1n-1}\in (x_{00})+Q$ then
   $u\in Q_1$.
\end{center}
Changing by elements of $Q_1$ and using the 
first claim, we may assume that $u$ is  
of the form $u=\sum_{b\in B} a_b b$ with almost all $a_b=0$.
By assumption we have $g(u)ts^{n-1}=g(ux_{1n-1})\in(zt^m)$,
hence $g(u)\in(zt^{m-1})$.
However, we have seen in the proof of the first claim
that this implies $a_b=0$ for all $b\in B$ and proves 
the second claim.

Finally, we prove our assertion about $\Hom_S(I,S)$:
Since $S=R/Q$ is a domain, the element $x_{00}$ is $S$-regular.
Moreover, since $I$ is an ideal of $S$ and $x_{00}$ is $S$-regular,
it follows that the $S$-module homomorphism
$\Hom_S(I,S)\to S$ given by $f\mapsto f(x_{00})$, is
injective.
Now, let $f\in\Hom_S(I,S)$ and set $u=f(x_{00})\in S$.
We are done if we show that $u$ lies inside the ideal generated
by $x_{00}$ and $x_{01}$ of $S$.
However, this follows from the computation
$$
ux_{1n-1}=f(x_{00})x_{1n-1}=f(x_{00}x_{1n-1})=x_{00}f(x_{1n-1})\in (x_{00})\subseteq S.
$$
together with the second claim above.
\qed

\begin{Remark}
  In fact, $\imath$ and $\phi$ are defined over the integers.
  Since they generate $\Hom_S(I,S)$ over any field, in particular
  over all prime fields, it follows
  that Proposition \ref{generators} holds in fact
  over the integers.
\end{Remark}
\medskip

Since $\Hom_S(I,S)$ has two generators in degree zero,
the Kustin--Miller unprojection with respect to the whole 
$S$-module $\Hom_S(I,S)$ yields a graded ring, whose component
of degree zero is a vector space of dimension at least
two, i.e., the unprojection ring is ''not geometric``.
Although even negatively graded rings occur in the description
of $3$-fold flips \cite[Section 11]{kinosaki},
we will use natural submodules of $\Hom_S(I,S)$ instead.

A geometric interpretation why the unprojection ring 
associated to the whole $S$-module $\Hom_S(I,S)$ 
does not give the ''right`` object is the following
observation:
the unprojection locus $\Gamma=V(I)$ is a curve with
self-intersection zero, whereas for the existence of a morphism
contracting $\Gamma$ we would need that $\Gamma$
has negative self-intersection.

\section{Generalised unprojections}
\label{sec:unprojection}

We keep the notations introduced so far.
As already noted before, taking the unprojection ring
with respect to the whole of $\Hom_S(I,S)$ yields
a graded ring, which is not ''geometric``,
which is why we consider suitable submodules.

In view of the natural short exact sequence  (\ref{poincare residue})
we will consider submodules of $\Hom_S(I,S)$ of the form
${\rm res}^{-1}(N)$, where $N$ is a submodule of
${\rm Ext}^1_S(S/I,S)$.
Recall that $\Hom_S(I,S)$ is generated as $S$-module by two elements
$\imath,\phi$ in our setup by Proposition \ref{generators}.
Then, in case $N$ is a cyclic $S$-module we are led to considering
submodules of $\Hom_S(I,S)$ that are generated by $\imath$
and another element $f\phi$,
where $f\in S$ is a homogeneous element.
Motivated by \cite{pr} and \cite{towards} we define

\begin{Definition}
 Let $S$ be the homogeneous coordinate ring of $\FF(m,n)$ inside
 $\PP^{m+n+1}$ and let $f\in S$ be homogeneous of
 degree $k\geq1$.
 The {\it generalised unprojection ring} of $S$ with respect to the 
 unprojection ideal $I$ and to $f$ is defined as
 $$
   S_{\rm un}(f) \,=\, \frac{S[T]}{(Tu-f\phi(u),\,u\in I)}\,,
 $$
 where $T$ is a variable of degree $k$.
\end{Definition}

\begin{Lemma}
 \label{avoidI}
 Let $f_1, f_2$ be homogeneous elements of $S$ of the same degree
 with $f_1-f_2\in I$.
 Then there exists an isomorphism of graded rings
 $$
    S_{\rm un}(f_1)\,\iso\,S_{\rm un}(f_2)\,.
 $$
 In particular, if $f_1\in I$ then $S_{\rm un}(f_1)$ is not 
 a domain.
\end{Lemma}

\proof
Let $i=f_2-f_1$, which is an element of $I$ by assumption.
Then for all $u\in I$ we calculate
$$
Tu-f_2\phi(u) \,=\, Tu-(f_1+i)\phi(u)\,=\,
\left(T-\phi(i)\right)u-f_1\phi(u).
$$
Thus a change of variables from $T$ to $T-\phi(i)$ 
(note that both elements are of the same degree) yields 
the desired isomorphism of graded rings.
\qed

\begin{Remark}
 That $S_{\rm un}(f)$ depends only on the submodule
 of $\Hom_S(I,S)$ generated by $\imath$ and $f\phi$
 and not on the particular choice of generators
 also follows from the intrinsic setup of \cite[Section 2]{towards}.
\end{Remark}

The previous lemma thus tells us that there is no loss
of generality choosing $f$ to be a homogeneous polynomial
in $x_{0m}$ and $x_{1n}$.
Since $x_{0m}$ and $x_{1n}$ are homogeneous coordinates
on $\Gamma=V(I)$ we remark the following.

\begin{Remark}
 \label{choiceofpoints}
 If the ground field is algebraically closed then 
 choosing a homogeneous element $f$ of degree $k$ 
 is equivalent to choosing $k$ points (counted with
 multiplicities) on $\Gamma$.
\end{Remark}

\begin{Theorem}
 \label{unprojection theorem}
 Let $f\neq0$ be homogeneous of degree $k\geq1$ in 
 $x_{0m}$ and $x_{1n}$. 
 Then $\Proj S_{\rm un}(f)$ is an integral, normal,
 and projectively Cohen--Macaulay surface inside
 weighted projective space $\PP(1^{m+n+2},k)$.
 Its homogeneous ideal is generated by
 the $2\times 2$-minors of the matrix
 \begin{equation}
 \label{unprojection equation}
 \left(
 \begin{array}{ccc|ccc|c}
   x_{00} & ... & x_{0m-1} & x_{10} & ... & x_{1n-1} & f\\
   x_{01} & ... & x_{0m} & x_{11} & ... & x_{1n} & T
 \end{array}
 \right)\,,
 \end{equation}
 where the $x_{ij}$ are of degree one and $T$ is of degree $k$.
\end{Theorem}

\proof
The description of the homogeneous ideal follows directly from the 
presentation of $S$ as the vanishing of the 
$2\times2$ minors of (\ref{scroll equation}) 
and the definition of $S_{\rm un}(f)$.

Let $Q_2$ be the ideal of $R_2=k[x_{ij},T]$ generated by the 
$2\times2$ minors of (\ref{unprojection equation}).
We want to show that if $P$ is a minimal prime ideal over $Q_2$ 
then it has codimension equal to $m+n$, which is
the maximum possible by a result of Eagon, compare
\cite[Exercise 10.9]{eis}.

Denote by $I^e$ the ideal of $R_2$ generated by the subset $I+Q_2$.
First, assume that $I^e\subseteq P$, which implies
${\rm codim}(P)\geq{\rm codim}(I^e)$.
However, $I^e$ contains $x_{00},...,x_{0m-1}$, 
$x_{10},...,x_{1n-1}$, $fx_{0m}$, 
which form a regular sequence, which implies
${\rm codim}(I^e)\geq m+n+1$.
Hence this case does not exist and we have
$I^e\not\subseteq P$, i.e.,
$V(P)\cap (V(Q_2)-V(I^e))\neq\emptyset$.
By \cite[Remark 2.5]{parallel}, the inclusion 
of rings induces an isomorphism 
$\Spec S_{\rm un}(f)-V(I^e)\iso\Spec S-V(I)$.
Since $\Spec S-V(I)$ is irreducible of dimension three, we
see that $P$ has codimension $m+n$.
Thus, every minimal prime over $Q_2$ 
has codimension $m+n$, which means that $S_{\rm un}(f)$ 
is a determinantal ring and such rings are known to be
Cohen--Macaulay by a result of Eagon, 
compare \cite[Theorem 18.18]{eis}.

By the above arguments, the irreducible open subset
$\Spec S_{\rm un}(f)-V(I^e)$ of $\Spec S_{\rm un}(f)$ meets every
irreducible component of $\Spec S_{\rm un}(f)$.
From this we conclude that $\Spec S_{\rm un}(f)$ is irreducible.

From the isomorphism $\Spec S_{\rm un}(f)-V(I^e)\iso\Spec S-V(I)$
it follows that $S_{un}(f)$ is generically reduced.
In particular, being Cohen--Macaulay there are no embedded components
and it follows that $S_{\rm un}(f)$ is reduced.
Together with the irreducibility it follows that $S_{\rm un}(f)$ is
a domain.

Using the isomorphism $\Spec S_{\rm un}(f)-V(I^e)\iso\Spec S-V(I)$
once more, we obtain normality outside $V(I^e)$.
A straightforward calculation using the Jacobian criterion 
shows normality along $V(I^e)$.
Thus, $S_{\rm un}(f)$ is normal.
\qed\medskip


\begin{Theorem}
 \label{main theorem}
  Let $\varpi:\tilde{X}\to X=\Proj S$ be the blow-up of the ideal
  $(I,f)$.
  Then there exists a factorisation
  $$\begin{xy}
   \xymatrix{
   &\tilde{X}
   \ar[dl]^{\varpi} \ar[dr]^{{\rm cont}_{\hat{\Gamma}}}\\
   X\,=\,\Proj S \ar@{-->}[rr] & & X_{\rm un}\,=\,\Proj S_{\rm un}(f)\,,
   }
  \end{xy}$$
  where ${\rm cont}_{\hat{\Gamma}}$ is contraction of the
  strict transform of $\Gamma=V(I)$ on $\tilde{X}$.
\end{Theorem}

\proof
Considered as an ideal of $R=k[x_{ij}]$, the
ideal $(I,f)$ is generated by the regular sequence
$x_{00},...,x_{0m-1},x_{10},...,x_{1n-1},f$.
By \cite[Exercise IV-26]{eh}, the Rees algebra
$\tilde{R}$ of $R$ with respect to $(I,f)$ is 
isomorphic to
$$
 R[T_{00},...,T_{0m-1},T_{10},...,T_{1n-1},T_f]/B,
$$
where the $T_{ij}$ and $T_f$ are indeterminants and 
where $B$ is the ideal generated by the $2\times2$ minors
of
$$
\left(
\begin{array}{ccc|ccc|c}
 x_{00}&...&x_{0m-1}&x_{10}&... &x_{1n-1}&f\\
 T_{00}&...&T_{0m-1}&T_{10}&... &T_{1n-1}&T_f
\end{array}
\right)\,.
$$
Taking $\Proj$, we obtain the blow-up $\varpi:\tilde{X}\to X$.
We denote by $\hat{\Gamma}$ the strict transform of $\Gamma$,
which is cut out by
$x_{00}=...=x_{0m-1}=0$, $x_{10}=...=x_{1n-1}=0$,
$T_{00}=...=T_{0m-1}=0$, and $T_{10}=...=T_{1n-1}=0$.

On the level of commutative algebra, $\varpi$ corresponds to 
eliminating the $T_{0i}$'s, the $T_{1j}$'s {\it and} $T_f$, 
whereas eliminating only the $T_{0i}$'s and the $T_{1j}$'s
induces a map ${\rm cont}_{\hat{\Gamma}}$ from $\tilde{X}$ onto
$X_{\rm un}$.

A straightforward calculation shows that ${\rm cont}_{\hat{\Gamma}}$
is in fact a morphism, that it is an isomorphism outside $\hat{\Gamma}$
and that it contracts $\hat{\Gamma}$ to the vertex of the weighted
projective space in which $X_{\rm un}$ lies.
Since $X_{\rm un}$ is normal by 
Theorem \ref{unprojection theorem},
Zariski's main theorem shows that
${\rm cont}_{\hat{\Gamma}}$ is in fact the contraction
of $\hat{\Gamma}$.
\qed\medskip

Let us assume that $D=\sum_i k_iP_i$ where the $k_i$ are positive
integers with $k=\sum_i k_i$ and where the $P_i$ are distinct points that
are rational over the ground field.
Note that this assumption on $D$ can always be fulfilled if the ground field 
is algebraically closed.
Then calculations similar to those in \cite[Chapter IV.2.3]{eh} show that 
we obtain $X_{\rm un}=\Proj S_{\rm un}(f)$ as follows:
\begin{enumerate}
 \item For each $i$, blow up $\Proj S$ at $P_i$. 
  Then blow up the intersection
  point of the strict transform of $\Gamma$ with the resulting $(-1)$-curve of
  the blow-up etc. until, for every $i$ we get a chain $C_i$ of $(k_i-1)$ 
  rational $(-2)$-curves and a $(-1)$-curve.
  It is understood that $C_i$ is empty if $k_i=1$.
 \item The strict transform $\hat{\Gamma}$ is a rational curve with 
  self-intersection $-k$.
  Contracting $\hat{\Gamma}$ and all the $C_i$'s we obtain $X_{\rm un}$.
\end{enumerate}
From this description we can read off the singularities of $X_{\rm un}$:
it has a toric singularity of type $\frac{1}{k}(1,1)$ coming from the 
contraction of $\hat{\Gamma}$, and every contracted chain $C_i$ contributes 
a cyclic quotient singularity of type 
$\frac{1}{k_i}(1,k_i-1)$, i.e., a Du~Val singularity of type $A_{k_i-1}$.

\section{Elementary transformations}

As in the previous sections, let $S$ be the homogeneous coordinate
ring of $\FF(m,n)$ inside $\PP^{m+n+1}$ given by (\ref{scroll equation})
and recall that we assumed $n\geq m\geq1$.
As unprojection divisor we take the line $\Gamma=V(I)$ lying on 
$\FF(m,n)$ as in (\ref{gamma equation}).
We note that $x_{0m}$ and $x_{1n}$ can be viewed as coordinates on $\Gamma$.
Moreover, in our unprojection setting we choose $0\neq f\in S$ of
{\em degree one}, which for our unprojection purposes we may assume
to be of the form
$$
f_{a,b}\,=\,a\,x_{0m}\,+\,b\,x_{1n}\,\mbox{ \qquad with \quad }[a:b]\,\in\,\PP^1,
$$
cf. Lemma \ref{avoidI}.
As already noted in Remark \ref{choiceofpoints}, our unprojection data
consists of an unprojection locus, which is a line, and a point on this 
line.

\begin{Proposition}
  There exists an isomorphism
  $$
    \Proj S_{\rm un}(f_{a,b})\,\iso\,\left\{
     \begin{array}{ll}
      \FF(m,n+1) & \mbox{ if } [a:b]=[0:1] \\
      \FF(m+1,n) & \mbox{ else,} 
     \end{array}\right.
  $$
  which is induced by a projective linear transformation of the ambient
  $\PP^{m+n+2}$.
\end{Proposition}

\proof
If $a=0$ or $b=0$ this follows directly from comparing
(\ref{unprojection equation}) with (\ref{scroll equation}).
We may thus assume $a\neq0$.
For all $i=0,...,n-1$
we add $b$ times the $(n-i)$.th column of the middle block of
(\ref{unprojection equation}) to the $(m-i)$.th column of the left block
of (\ref{unprojection equation}).
This is possible since we assumed $n\geq m$ and a linear change
of variables yields the desired isomorphism.
\qed

\begin{Remark} 
  By Theorem \ref{main theorem}, we have realised all elementary transformations 
  of Hirzebruch surfaces in our setting.
\end{Remark}

Moreover, we obtain a family of unprojections parametrised
by $\PP^1$.
One member of this family is isomorphic to $\FF_{n-m+1}$, whereas
all the others are isomorphic to $\FF_{n-m-1}$.
This fits nicely into the deformation and degeneration 
theory of Hirzebruch surfaces as explained in
\cite[Theorem VI.8.1]{bhpv}.

The inverse of the unprojection 
$\Proj S\dashrightarrow\Proj S_{\rm un}(f_{a,b})$, corresponds to 
eliminating the new variable $T$ of $S_{\rm un}(f_{a,b})$, 
and is induced by a projection from $\PP^{m+n+2}$ onto $\PP^{m+n+1}$,
confer \cite[Proposition 8.20]{har}.

\section{(Bi)canonical images of Horikawa surfaces}

In this final section we work over the complex numbers.
If $S$ is a minimal surface of general type then 
Noether's inequality $K^2\geq2p_g-4$ holds true,
confer \cite[Theorem VII.3.1]{bhpv}.
In case of equality $K^2=2p_g-4$, i.e., if $S$ is a
so-called {\it even Horikawa surface}, then
the canonical map
is a generically finite morphism of degree $2$ onto a surface
of minimal degree in $\PP^{p_g-1}$, which is the key to the 
classification of these surfaces,
compare \cite[Chapter VII.9]{bhpv}.
Also, it is not difficult to show that the bicanonical map
is a morphism that coincides with the canonical map
followed by the second Veronese embedding.

Now, let $S$ be an {\it odd Horikawa surface}, i.e., a minimal
surface of general type with $K^2=2p_g-3$.
These have been classified in \cite{horikawa}, and we will
assume that we are in case A in Horikawa's terminology
with smooth canonical image.
This is the generic case, and if $p_g\geq7$ it is even 
automatically fulfilled, compare \cite[Theorem 1.3]{horikawa}.

Then the image of the canonical map is a smooth surface $X$ 
of minimal degree in $\PP^{p_g-1}$.
Hence there exist integers $n\geq m\geq 1$ with
$p_g=m+n+2$ such that the canonical image is $\FF(m,n)$, which 
is abstractly isomorphic to $\FF_{n-m}$.
The canonical system $|K_S|$ has a unique base point, whose
indeterminacy is resolved by a single blow-up $\tilde{S}\to S$.
The resulting $(-1)$-curve on $\tilde{S}$ maps to a line
$\Gamma\subset\FF(m,n)\subset\PP^{p_g-1}$ and after an appropriate choice of
coordinates we may assume that $\Gamma$ and $\FF(m,n)$
are as in Section \ref{sec:scrolls}.
Then $S$ determines two  points $x,y$ on $\Gamma$, possibly infinitely 
near, cf. \cite[Theorem 1.3]{horikawa}. 
We denote by $\pi:\tilde{X}\to X$ their blow-up,
by $E_x$, $E_y$ the corresponding exceptional divisors, and
by $\hat{\Gamma}$ the strict transform of $\Gamma$ on $\tilde{X}$.
Moreover, we obtain a factorisation $\tilde{S}\to S^\ast\to \tilde{X}\to X$,
where $\tilde{S}\to S^\ast$ is birational, $S^\ast$ has at
worst Du~Val singularities, and where $S^\ast\to\tilde{X}$
is finite and flat of degree $2$.
On $\tilde{X}$ we consider the line bundle
$$
{\cal L}\,=\,  \OO_{\tilde{X}}(\pi^\ast\Delta_0\,+\,(n-4)\pi^\ast\Gamma\,-\,2E_x\,-\,2E_y)\,.
$$
Then the canonical map of $\tilde{S}$ (and hence $S$) factors over
the complete linear system $|{\cal L}|$ on $\tilde{X}$ and we
already noted that we can identify its image with $\FF(m,n)$.
Moreover, the $(-1)$-curve on $\tilde{S}$ maps to $\hat{\Gamma}$
on $\tilde{X}$ and thus maps to $\Gamma$ under the canonical map.

From \cite[Theorem 1.3]{horikawa} it follows easily that the bicanonical
map factors over $|{\cal L}^{\otimes 2}|$ on $\tilde{X}$.
This map contracts $\hat{\Gamma}$ to an $A_1$-singularity.
A straightforward computation counting the quadratic
relations coming the $2\times2$ minors of (\ref{scroll equation}),
we see that the map $H^0({\cal L})^{\otimes 2}\to H^0({\cal L}^{\otimes2})$
has $1$-dimensional coimage.
This implies that the canonical ring of $S$ has $p_g$ generators
in degree one and precisely one new generator in degree two.

Let us denote by $X_1\subset\PP(1^{p_g})$ and by 
by $X_2\subset\PP(1^{p_g},2)$
the projection from the canonical model of $S$ onto the weighted
projective space corresponding to generators in degree $1$ (canonical
image) and generators in degree $\leq2$ (weighted bicanonical image).
Then, putting all observations above together we can interpret 
Horikawa's results \cite[Section 1]{horikawa} in our setting as
follows:

\begin{Proposition}
 The odd Horikawa surface $S$ determines a line $\Gamma$ on $X_1$
 and two points $\{x,y\}$ on this line, which are possibly infinitely near.
 Then the inverse of the natural projection $X_2\to X_1$,
 $$
     X_1\,\dashrightarrow\,X_2
 $$
 is a generalised unprojection with unprojection data
 $\Gamma$ and divisor $D=x+y$ on $\Gamma$.
\end{Proposition}
\medskip

Finally, we remark that we cannot realise the other birational
modifications of Theorem \ref{main theorem} 
by images of pluricanonical maps of surfaces of general type:
If the canonical map has two-dimensional image then
$p_g\geq3$ and in order to get a scroll as canonical image
we even need $p_g\geq4$.
For such surfaces we have $K^2\geq4$ by 
Noether's inequality.
However, for minimal surfaces of general type $X$ with $K_X^2\geq3$
all pluricanonical images ${\rm im}(\varphi_i(X))$ with $i\geq3$ 
are birational to $X$ by
Bombieri's theorem (confer \cite[Theorem VII.5.1]{bhpv})
and thus no rational surfaces.

\bigskip
\noindent
Christian Liedtke \\
Department of Mathematics\\
Stanford University\\
450 Serra Mall\\
Stanford, CA 94305, USA\\
e-mail: liedtke@math.stanford.edu

\bigskip
\noindent
Stavros Papadakis, \\
Center for Mathematical Analysis, Geometry, and Dynamical Systems\\
Departamento de Matem\'atica, Instituto Superior T\'ecnico\\
Universidade T\'ecnica de Lisboa\\
1049-001 Lisboa, Portugal\\
e-mail: papadak@math.ist.utl.pt

\end{document}